\documentclass[a4paper,leqno,11pt]{amsart}

\usepackage[T1]{fontenc}
\usepackage[sc,osf]{mathpazo}
\usepackage{setspace}
\usepackage{csquotes}
\usepackage{amsmath,amsthm}
\usepackage{amssymb}
\usepackage{xspace}
\usepackage{bm}
\usepackage{amstext}
\usepackage{comment}
\usepackage{amsfonts}
\usepackage{graphicx}
\usepackage[mathscr]{euscript}
\usepackage{amscd}
\usepackage{latexsym}
\usepackage{amssymb}
\usepackage{enumerate}
\usepackage{mathrsfs}
\usepackage[cmtip,all]{xy}
\usepackage{xcolor}
\usepackage{standalone}
\usepackage{tikz}
\usetikzlibrary{arrows.meta}
\usepackage{tikz-cd}
\usetikzlibrary{intersections}
\usetikzlibrary{calc}
\usetikzlibrary{decorations.markings}
\usetikzlibrary{arrows}
\usetikzlibrary{shapes,arrows,shapes.multipart}

\usepackage{xcolor}

\usepackage{enumitem}

\usepackage[cmtip,all]{xy}

\newtheorem{thm}{Theorem}[section]
\newtheorem{prop}[thm]{Proposition}
\newtheorem{lemma}[thm]{Lemma}
\newtheorem{cor}[thm]{Corollary}

\theoremstyle{definition}
\newtheorem{defn}[thm]{Definition}

\newtheorem{remark}[thm]{Remark}

\newtheorem*{thma}{Theorem A}

\newtheorem*{thmb}{Theorem B}

\newcommand{\uM}{\underline{M}}
\newcommand{\uN}{\underline{N}}
\newcommand{\bF}{\mathbb{F}_2}

\newcommand{\ug}{\underline{g}}

\newcommand{\upi}{\underline{\pi}}
\newcommand{\Sp}{\mathsf{Sp}}

\newcommand{\res}{\mathrm{res}}

\newcommand{\F}{\bF}

\newcommand{\tH}{\tilde{H}}
\newcommand{\cP}{\mathcal{P}}

\newenvironment{mysubsection}[1]
{%
  \refstepcounter{thm}%
  \par\medskip
  \noindent\textbf{\thethm.\ #1.}\upshape\ignorespaces
}
{\par\medskip}

\newenvironment{myeq}
{%
  \refstepcounter{thm}%
  \begin{equation}
  \tag{\thethm}
}
{%
  \end{equation}
}

\newcommand{\rb}{\bar\rho}

\newcommand{\coker}{\operatorname{coker}}

\begin{document}
\title{Characterizing $2$-slices over finite groups}

\author{Surojit Ghosh}
\address{Department of Mathematics, Indian Institute of Technology, Roorkee, Uttarakhand-247667, India}
\email{surojit.ghosh@ma.iitr.ac.in; surojitghosh89@gmail.com}

\author{Ankit Kumar}
\address{Department of Mathematics, Indian Institute of Technology, Roorkee, Uttarakhand-247667, India}
\email{ankit\_k@ma.iitr.ac.in}

\date{\today}
\subjclass{55N91, 55P91 (primary); 57S17, 55Q91 (secondary)}
\keywords{Equivariant homotopy theory, Slice tower, Mackey functor}
\date{\today}
\subjclass{55N91, 55P91 (primary); 57S17, 55Q91 (secondary)}
\keywords{Equivariant homotopy theory, Slice tower, Mackey functor}

\begin{abstract}
For a finite group $G$, we characterize the $2$-slices of $G$-spectra in
the regular slice filtration in terms of their homotopy Mackey functors
and $k$-invariant. We explicitly compute this for cyclic $p$-groups and the alternating group $A_4$.
\end{abstract}

\maketitle

\section{introduction}

The slice filtration generalizes Dugger’s filtration \cite{Dug05} and is modelled on Voevodsky’s motivic slice filtration \cite{Voe02}.  This filtration is an
equivariant analogue of the Postnikov filtration. It replaces the spheres
$S^n$ by the slice cells $G_+\wedge_HS^{k\rho_H}$, indexed by representation-theoretic dimension, and decomposes a $G$-spectrum $X$
into slices $P^n_nX$. The resulting slice spectral sequence was the main
computational tool in their solution of the Kervaire invariant one
problem \cite{HHR16}, and it has since been used to compute many equivariant homotopy
groups. However, all the calculations here are for the regular slices as introduced by Ullman \cite{U}. Each such computation starts by identifying the slices, so it is
useful to know which $G$-spectra are $n$-slices.

For small $n$ this is known. A $G$-spectrum is a $0$-slice if and only if
it is an Eilenberg--Mac~Lane spectrum $H\uM$, and $\Sigma^1H\uM$ is a
$1$-slice if and only if all restrictions of $\uM$ are injective
\cite{HHR16}. Ullman \cite{U} characterized the $(-n)$-slices of the form
$\Sigma^{-1}H\uM$. Since $\Sigma^{\rho_G}$ shifts slices by $|G|$, these
results describe all slices in degrees $0, \cdots, |G|-1\pmod{|G|}$. For the Klein four-group, the open case is that of $2$-slices. Slone settled it in \cite{Slo22}.

In this note we characterize the $2$-slices for every finite group $G$.
For $H\le G$ and a $G$-Mackey functor $\uM$, write
$\varphi_H\colon\uM(G/H)\to\tH^0_H(S(\bar\rho_H)_+;\uM)$ for the map
induced by the collapse map $S(\bar\rho_H)\to\mathrm{pt}$, where
$\bar\rho_H$ is the reduced regular representation.

\begin{thma}\label{mainA}
Let $G$ be a finite group.
\begin{enumerate}[label=\textup{(\alph*)}]
\item $\Sigma^2H\uM$ is a $2$-slice if and only if $\res^H_e$ is injective
for all $|H|=2$ and $\varphi_H$ is bijective for all $|H|\ge3$.
\item $\Sigma^1H\uM$ is a $2$-slice if and only if, for every $H\le G$,
the map $\uM(G/H)\to\displaystyle\prod_{C\le H,\,|C|=2}\uM(G/C)$ is injective. In
particular, $\uM(G/H)=0$ whenever $|H|$ is odd.
\item A $G$-spectrum $Y$ is a $2$-slice if and only if $\upi_kY=0$ for
$k\ne1,2$, $\upi_1Y(G/e)=0$, $\Sigma^2H\upi_2Y$ is a $2$-slice, and for
every $H$ with $|H|\ge3$ the $k$-invariant of $Y$ induces an injection
\[
\bigcap_{L\lneq H}\ker\res^H_L\hookrightarrow
\tH^1_H\bigl(S(\bar\rho_H)_+;\upi_2Y\bigr).
\]
\end{enumerate}
\end{thma}

This theorem provides a practical criterion for determining whether a given spectrum is a $2$-slice. In particular, when the spectrum is of the form $\Sigma^2H\uM$ or $\Sigma^1 H\uM$, the criterion yields explicit necessary and sufficient conditions for it to be a $2$-slice. We investigate these cases in Section~4. 


To make Theorem~A applicable, we compute the groups appearing in its statement. The cell structure of $S(\bar\rho_G)$, as described in \cite{HHR16}, leads to the following formula.

\begin{thmb}\label{mainB}
$\tH^0_G\bigl(S(\bar\rho_G)_+;\uM\bigr)$ is the group of families
$(x_L)_{L\lneq G}$, with $x_L\in\uM(G/L)$, that are compatible with
restrictions, and with conjugations at subgroups of index at least $3$.
\end{thmb}
See Theorem~\ref{zerosrho}.

\medskip

 Note that $\tH^0_G(S(\bar\rho_G)_+;\uM)$ can be strictly larger than
$\tH^0_G(E\mathcal P_+;\uM)=\lim_{\mathcal O_{\mathcal P}(G)^{op}}\uM$ \cite{MNN19}.
We also express Theorem~B in terms of maximal subgroups and their
Weyl groups. From this we obtain explicit descriptions of the $2$-slices.

\begin{mysubsection}{Organization} ~ We organize the remainder of the paper as follows. In Section~3, we present a combinatorial description of $0^{th}$-Bredon cohomology of the sphere in the reduced regular representation of $G.$ In Section~4, we provide a general classification result for $2$-slices and consequently we describe this result for the cyclic $p$-groups and the alternating group $A_4.$
\end{mysubsection}

\paragraph{\textbf{Acknowledgement.}}~The research of Surojit Ghosh was supported by ANRF, India, with grant number  ARG/2025/000434/MS. Ankit Kumar's research is supported by the University Grants Commission, India.

\section{Foundation}
In this section, we briefly review slice filtration and slice tower for a $G$-spectrum $X$. Throughout, we use the regular slices introduced by Ullman in \cite{U}. Let $\Sp^G$ denote the category of genuine $G$-equivariant spectra. A
\emph{slice cell} is a $G$-spectrum of the form $G_+\wedge_H S^{k\rho_H}$ with
$H\le G$ and $k\in\mathbb{Z}$; its \emph{dimension} is $k|H|$.   Here $\rho_H$ denotes the regular representation of the subgroup $H \le G.$ For each integer $n$, define $\tau_{\ge n} \subseteq \Sp^G$ to be the localizing subcategory generated by the slice cells.  By definition, $\tau_{\ge n}$ is closed under taking weak equivalences, cofibers, extensions and arbitrary wedges. 

\begin{defn}
    We write $X \ge n$ if $X \in \tau_{\ge n},$ and $X\le n$ if $[C, X]^G=0$ for all $C \in \tau_{n+1}.$ Equivalently, it gives $[G_+\wedge_H S^{k\rho_H+r}, X]^G=0$ for all $H \le G$, $r\ge 0$, and $k|H|>n.$ 
\end{defn} 

 Associated to $\tau_{\ge n}$ there is a natural localization functor $P^{n-1}$ on $ \Sp^G$ whose acyclic objects are exactly the objects of $\tau_{\ge n}.$ Since $\tau_{\ge n+1} \subseteq  \tau_{n}$, there is a natural transformation $P^n \to P^{n-1}$, which yields the slice tower $\{P^n X\}_{n\in \mathbb{Z}}$ of each $G$-spectrum $X.$  The homotopy fiber of the map $P^n X \to P^{n-1} X$, denoted by $P^n_nX$, is called the \emph{$n$-slice} of $X.$ Equivalently, a spectrum $X$ is an $n$-slice if $X\le n$ and $X \ge n.$ 

 Hill and Yarnall provided the algebraic characterization of $X \ge n$ in terms of the geometric fixed points. We write $\Phi^HX$ for the geometric $H$-fixed point of the restriction $\res^G_HX \in \Sp^H$, which is a symmetric monoidal functor of $\Sp^H.$

\begin{thm}[\cite{HY18}]{\label{lbound}}  
A $G$-spectrum $X$ satisfies $X \ge n$ if and only if $\pi_k\Phi^HX=0$ for all $H\le G$ and $k < n/|H|.$ 
\end{thm}
In particular, one observes $\Sigma^mH\underline M\ge m$ for $m\ge0$ since the geometric fixed point of $\Sigma^m H\uM$ is $\Sigma^m\Phi^LH\uM$, hence $(m-1)$-connected and $m/|L|\le m.$

The slices in low dimensions are described in terms of Mackey functors.

\begin{prop}[\cite{HHR16,U}]{\label{oldchar}}Let $\uM$ be a $G$-Mackey functor and $X \in \Sp^G.$
\begin{enumerate}
    \item $X$ is a $0$-slice if and only if $X\simeq H\underline M$ for some $G$-Mackey functor $\uM.$
    
    \item $\Sigma^1H\underline M$ is a $1$-slice if and only if all restrictions of $\underline M$ are injective.
    
    \item $\Sigma^{-1}H\underline M$ is a $(-n)$-slice if and only if transfers out of subgroups of order $\ge n$ are surjective and $\underline M(G/H)=0$ for $|H|<n$.
\end{enumerate}
\end{prop}

Another exciting property that Hill, Hopkins, and Ravenel observed is that the slices are $\rho_G$-periodic. More precisely,
\begin{prop}[{\cite[Corollary~4.25]{HHR16}}]
 If \(X\) is an \(n\)-slice, then \(\Sigma^\rho X\) is an \((n+|G|)\)-slice. i.e.
\[
P^{k+|G|}_{k+|G|}(\Sigma^{\rho_G}X)\simeq\Sigma^{\rho_G}P^k_k(X)
\]
where $\rho_G$ is the regular representation of group $G$.
\end{prop}

\section{limit presentation of $\tilde{H}_G^0(S({\bar{\rho}_G}); \uM)$}
In this section, we prove certain presentations of the $0^{th}$-Bredon cohomology of the sphere in the reduced regular representation with coefficients in a Mackey functor $\uM$. 

To begin with, let us consider a Mackey functor $\uM$ and  the cofiber sequence
\[
S(\bar{\rho}_G)_+ \to S^0 \to S^{\bar{\rho}_G}.
\]
This yields the associated long exact sequence
\[
0 \to \tilde{H}_G^0(S^{\bar{\rho}_G}; \uM) \to \uM(G/G) \stackrel{\varphi_G}{\to}  \tilde{H}_G^0(S(\bar{\rho}_G)_+; \uM) \to \tilde{H}_G^1(S^{\bar{\rho}_G}; \uM) \to 0
\]
in Bredon cohomology with $\uM$-coefficients. Then  by \cite[Corollary 3.6]{HHR16}, the kernel is
\begin{myeq}\label{kervarphi}
\ker(\varphi_G (\uM))\cong \tilde{H}_G^0(S^{\bar{\rho}_G}; \uM) \cong \bigcap_{L \in \mathcal{P}_G} \ker (\res^G_{L}(\uM))
\end{myeq}
where $\cP_G$ is the poset of all non-trivial proper subgroups of $G.$

\begin{lemma}\label{krho}
For $k\ge 1$, there is an isomorphism $\tilde{H}_G^0(S^{k\bar{\rho}_G}; \uM) \cong \ker(\varphi_G)$. 
\end{lemma}

\begin{proof}
    We will prove this by induction on $k.$ Let us consider the cofiber sequence
    \[
    S^{m\bar{\rho}_G}\hookrightarrow S^{(m+1)\bar{\rho}_G}\to C
    \]
    Denote by $d_L= \dim \bar{\rho}_G^L$ for $L \in \cP_G$. Then $C^L\simeq S^{(m+1)d_L} \vee S^{md_L+1}$ and $C^G\simeq \ast.$ Hence the result follows.
\end{proof}

\begin{thm}\label{zerosrho}
Let $G$ be a finite group and $\uM$ a $G$-Mackey functor. Then
$\tH^0_G\bigl(S(\bar{\rho}_G)_+;\uM\bigr)$ is isomorphic to the group of
families
\[
(x_L)_{L\in \mathcal{P}_G}\in\prod_{L\in \mathcal{P}_G}\uM(G/L)
\]
satisfying the following two conditions:
\begin{enumerate}[label=\textup{(\roman*)}]
\item $\res^{L_2}_{L_1}x_{L_2}=x_{L_1}$ for all $L_1\le L_2\lneq G$;
\item $c_g\,x_L=x_{gLg^{-1}}$ for all $g\in G$ and all $L\in \mathcal{P}_G$ with
$|G:L|\ge 3$.
\end{enumerate}
\end{thm}

\begin{proof}
   Note that 
$$\tH^0_G(X_+;\uM)=H^0_G(X;\uM)\cong\mathrm{Hom}\big(\mathbb Z[\pi_0X^{(-)}],\uM\big).$$ 
Thus a class is a family of functions
$\eta_L\colon\pi_0X^L\to\uM(G/L)$ with
\begin{myeq}\label{erel}
\eta_{L_1}(c')=\res^{L_2}_{L_1}\eta_{L_2}(c),\qquad
\eta_{gLg^{-1}}(gc)=c_g\,\eta_L(c),
\end{myeq}
where $L_1\le L_2$ and $c'\in\pi_0X^{L_1}$ is the component containing
$c\in\pi_0X^{L_2}$.

For $X =S(\bar{\rho}_G),$ we have
\[
X^L\simeq \begin{cases}
    S^{[G:L]-2} & L \in \cP_G,\\
    \emptyset & \text{otherwise.} 
\end{cases}
\]
Hence $X^L$ is connected if and only if $[G:L]\ge3$; this holds whenever $L$ is properly contained in a proper subgroup. If $[G:H]=2$, then $H$ is normal and $X^H=\{\pm u_H\}$, with $g\notin H$ swapping the two points.

If $[G: L]\ge 3,$ then $\mathbb{Z}[\pi_0 X^L]\cong \mathbb{Z}.$ Let us call the generator by $w_L.$ Then define $x_L:= \eta_L(w_L)$. If $[G: L]=2$, then $\mathbb{Z}[\pi_0 X^L]\cong \mathbb{Z}\{u_L\}\oplus \mathbb{Z}\{-u_L\}$ with $g \in G\setminus L,$ such that $gu_L=-u_L.$ In this case, define $x_L:=\eta_L(u_L)$.

Given a family of functions $\eta_L\colon \mathbb{Z}[\pi_0X^L]\to\uM(G/L)$. Define $\eta_L(-u_L)= gx_L$ for $g\notin L$, which is independent of $g$
since $c_h=\mathrm{id}$ on $\uM(G/L)$ for $h\in L$.
\end{proof}

 Let \(K_1,\ldots, K_r\) be representatives of the conjugacy classes of maximal proper subgroups of \(G\), and set \(W_i=N_G(K_i)/K_i\). By maximality, \(W_i=1\) if \(K_i\) is not normal in \(G\), whereas \(W_i=G/K_i\) is cyclic of prime order if \(K_i\) is normal. We now describe the component corresponding to each \(K_i\) in \(\tH^0_G\bigl(S(\bar{\rho}_G)_+;\uM\bigr)\). Define

\[
A_i=
\begin{cases}
\uM(G/K_i)^{W_i} & \text{if } [G:K_i]\ge 3,\\[3pt]
\bigl\{x\in\uM(G/K_i):\res^{K_i}_L(c_gx-x)=0
\text{ for all maximal } L<K_i\bigr\} & \text{if } [G:K_i]=2,
\end{cases}
\]
where $g\in G\setminus K_i$ in the second case.

\begin{cor}\label{zerosrho-max}
Let $G$ be a finite group and $\uM$ a $G$-Mackey functor. Then
\[
\tH^0_G\bigl(S(\bar{\rho}_G)_+;\uM\bigr)\;\cong\;
\Bigl\{(x_i)\in\prod_{i=1}^rA_i\ :\
\res^{K_i}_{K_i\cap gK_jg^{-1}}x_i
=\res^{gK_jg^{-1}}_{K_i\cap gK_jg^{-1}}c_gx_j
\ \text{ whenever } K_i\ne gK_jg^{-1}\Bigr\}.
\]
It suffices to take $g$ in a set of representatives of
$K_i\backslash G/K_j$. In particular, if $G$ has no subgroup of index
$2$, then $A_i=\uM(G/K_i)^{W_i}$ for all $i$.
\end{cor}
\begin{proof}  We begin by considering conditions~(i) and~(ii) of Theorem~\ref{zerosrho}. By condition~(i), it suffices to consider the components corresponding to maximal subgroups of \(G\). Furthermore, condition~(ii) allows us to restrict our attention to representatives of their conjugacy classes. Consequently, the compatibility conditions in Theorem~\ref{zerosrho} yield
\[
\res^{K_i}_{K_i\cap gK_jg^{-1}}x_i
=
\res^{gK_jg^{-1}}_{K_i\cap gK_jg^{-1}}c_gx_j.
\]
Hence the result follows.
\end{proof}
\begin{remark}{\label{zero-rem}}
 If every maximal subgroup of \(G\) is normal, then the equation in the above corollary reduces to condition~(1) of Theorem~\ref{zerosrho}.
\end{remark}

For the $p$-cyclic group $C_{p^n}$, the above corollary specialises to the following result.
\begin{cor}\label{algcpn}
    Let $G=C_{p^n}$ with $n\ge1$. Then
\[
\tH^0_G\bigl(S(\bar{\rho}_G)_+;\uM\bigr)\cong
\begin{cases}
\uM(G/C_{p^{n-1}})^{G/C_{p^{n-1}}} & p \text{ odd},\\[2pt]
\uM(G/e) & p=2,\ n=1,\\[2pt]
\bigl\{x\in\uM(G/C_{p^{n-1}}):\res^{C_{p^{n-1}}}_{C_{p^{n-2}}}(x)\in\uM(G/C_{p^{n-2}})^{G/C_{p^{n-2}}}\bigr\} & p=2,\ n\ge2,
\end{cases}
\]
and under this isomorphism $\varphi_G(m)=\res^G_{{G/C_{p^{n-1}}}}(m)$.
\end{cor}

\begin{proof}
Since ${C_{p^{n-1}}}$ is the unique maximal subgroup of $G$, and it is normal
of index $p$ with the Weyl group ${G/C_{p^{n-1}}}$. So by Corollary~\ref{zerosrho-max} the product $\prod_{i=1}^rA_i$ is $\uM(G/C_{p^{n-1}})^{G/C_{p^{n-1}}}$ if $p$ is odd.

Consider the case when  $p=2$. If $n=1$, then the result follows from the identification $S(\bar{\rho_{C_2}})= C_2/e$. If $n\ge2$, then $C_{p^{n-2}}$ is the unique
maximal subgroup of $C_{p^{n-1}}$, so $A_1$ is defined by
$c_g\res^{C_{p^{n-1}}}_{C_{p^{n-2}}}(x)=\res^{C_{p^{n-1}}}_{C_{p^{n-2}}}(x)$ for a generator $g$ of $G$. Hence the result follows.
\end{proof}

For the Klein four-group $K_4$, let $H_1$, $H_2$, and $H_3$ denote its three maximal subgroups, each of index two. Then Corollary~\ref{zerosrho-max} reduces to the following result.

\begin{cor}\label{algk4}
    Let $G=K_4$ and $\alpha_i$ be the non-trivial irreducible characters of $K_4$ with $\ker \alpha_i=H_i$ for $i=1,2,3.$ Then  
    \[
\tH^0_{K_4}\bigl(S(\rb_{K_4})_+;\uM\bigr)\;\cong\;
\bigl\{(x_1,x_2,x_3)\in \bigoplus_{i=1}^3\uM(K_4/H_i) :\
\res^{H_1}_ex_{1}=\res^{H_2}_ex_{2}=\res^{H_3}_ex_{3}\bigr\}.
\]
In particular, no $K_4/H_i$-invariance is imposed on $x_i$. Under this isomorphism, 
\[
\varphi_{K_4}(m)=(\res^{K_4}_{H_1}m,\res^{K_4}_{H_2}m,\res^{K_4}_{H_3}m).
\]
\end{cor}
\begin{proof}
The maximal proper subgroups of \(K_4\), up to conjugacy, are \(H_1,H_2,\) and \(H_3\). Hence, Corollary~\ref{zerosrho-max}, together with Remark~\ref{zero-rem}, gives the stated description of \(\tH^0_{K_4}\bigl(S(\rb_{K_4})_+;\uM\bigr)\). The description of \(\varphi\) follows directly from \eqref{kervarphi}.
\end{proof}
Next, we consider the alternating group \(A_4\), which has order \(12\). Since \(A_4\) has no subgroup of index \(2\), the fixed-point space \(S(\rb_{A_4})^J\simeq S^{[A_4:J]-2}\) is connected for every proper subgroup \(J<A_4\). Recall that the maximal proper subgroups of \(A_4\), up to conjugacy, are \(K_4\) and \(C_3\). Thus, Corollary~\ref{zerosrho-max} yields the following description of \(\tH^0_{A_4}\bigl(S(\rb_{A_4})_+;\uM\bigr)\).
\begin{cor}\label{alga4}
    For any $A_4$-Mackey functor $\uM$, we have
    \[
\tH^0_{A_4}\bigl(S(\rb_{A_4})_+;\uM\bigr)\;\cong\;
\bigl\{(x,y)\in \uM(A_4/K_4)^{C_3} \times \uM(A_4/C_3) : \res^{K_4}_ex=\res^{C_3}_ey\bigr\}.
    \]
    and the map $\varphi_G: \uM(A_4/A_4)\to \tH^0_{A_4}\bigl(S(\rb_{A_4})_+;\uM\bigr)$ is given by $(\res^{A_4}_{K_4}, \res^{A_4}_{C_3})$.
\end{cor}

\section{Characterization of $2$-slices}
The primary objective of this section is to establish a characterization of \(2\)-slices for arbitrary \(G\)-spectra. We subsequently apply this characterization to the alternating group \(A_4\) and \(p\)-cyclic groups.  We start with some classification results for the Eilenberg-Mac Lane spectra.
\begin{prop}\label{twosus}
   Let $G$ be a finite group and $\uM$ be a $G$-Mackey functor. The spectrum $\Sigma^2 H\uM$ is a $2$-slice if and only if the following conditions hold for each subgroup $H$.
   \begin{enumerate}
       \item If $|H|=2$, $\res^H_e$ is injective.
       \item If $|H|\ge 3$, then $\varphi_H:\res_H\uM(H/H) \to \tH^0_H(S(\bar{\rho}_H)_+; \res_H\uM)$ is bijective.
   \end{enumerate}
   \end{prop}

   \begin{proof} We first assume that conditions (1) and (2) hold.
   
      By Theorem~\ref{lbound}, \(\Sigma^m H\underline{M}\geq m\) for every \(m\geq 0\). In particular, we are interested in the case \(m=2\). To determine the conditions on \(\underline{M}\) under which \(\Sigma^2 H\underline{M}\) is a \(2\)-slice, we consider the groups
\[
[S^{k\rho_H+r},\Sigma^2 H\underline{M}]^H
\cong
\widetilde{H}^{2-k-r}_H(S^{k\bar{\rho}_H};\underline{M}),
\]
for \(k|H|>2\) and \(r\geq 0\). For $\Sigma^{2}H\uM\leq 2$, these groups must vanish. For this, note that if $2-k-r<0$, then this is evident. But when $ 2-k-r\ge 0$, we analyse this case by case.

\noindent Suppose that \(|H|=2\). The condition \(k|H|>2\) implies that \(k>1\). Thus, the only case in which \(2-k-r\geq 0\) is \((k,r)=(2,0)\). In this case, we have $\widetilde{H}^{0}_H(S^{2\bar{\rho}_H};\underline{M})
\cong \ker(\varphi_H)
=\ker(\res^H_e)$ from Lemma~\ref{krho}, which vanishes by (1) hypothesis.     

\noindent If $|H|\ge 3,$ then the possibilities for the pair $(k,r)$ are $(1,0), (1,1)$ and $(2,0).$ For the first case the group is $\tH^1_H(S^{\bar{\rho}_H}; \uM)\cong\coker \varphi_H.$ For the other cases, using Lemma \ref{krho}, it is again $\ker(\varphi_H)$. Hence both results follow by (2).

For the other direction, suppose that $\Sigma^2 H\underline{M}$ is a $2$-slice. 
Then conditions (1) and (2) follow directly from the definition of a $2$-slice.
   \end{proof}

   The proposition above provides a complete characterization of when 
$\Sigma^2 H\underline{M}$ is a $2$-slice. We next investigate an 
analogous algebraic criterion for $\Sigma^1 H\underline{M}$ to be 
a $2$-slice.

   \begin{prop}\label{onesus}
   Let $G$ be a finite group and $\uM$ be a $G$-Mackey functor. The spectrum $\Sigma^1 H\uM$ is a $2$-slice if and only if $\uM(G/e)=0$ and $\varphi_H$ is injective for  for every $H \in G$ with $|H|\ge 3.$
  \end{prop}

  \begin{proof} We first assume that the above conditions hold and show that
$\Sigma^1 H\underline{M}$ is a $2$-slice.
  
  Again for $\Sigma^{1} H\uM\leq2 $, we need the group $\tH^{1-k-r}_H(S^{k\bar{\rho}_H}; \uM)$ to be zero whenever $k |H|> 2$ and $r \ge 0.$ For this the only non-trivial pair $(k,r)$ is $(1,0)$, which implies the value of the group is $\ker \varphi_H$ for $|H|\geq 3$. Hence $\Sigma^{1} H\uM \le 2$ if and only if $\ker \varphi_H=0$ for all $|H|\ge 3.$

To establish that \(\Sigma^1 H\uM\geq 2\), we appeal to Theorem~\ref{lbound}, which implies that
\[
\pi_k(\Phi^e\Sigma^1 H\uM)=0
\qquad\text{for all }k<2.
\]
But \(\uM(G/e)=0\), which gives the desired result.

The other direction follows directly from the definition of a $1$-slice.
\end{proof}

\begin{remark}
From Proposition~\ref{oldchar}, $\Sigma^{1}H\uM$ is a
$1$-slice whenever all the restriction maps of $\uM$ are injective.
However, if $M(G/e)=0$, then $\Sigma^{1}H\uM$ is contractible.    
\end{remark}

The following proposition is an evident fact which follows from the definition. Slone \cite{Slo22} used this for $K_4$ without a proof; we generalise it for any finite group as follows.

\begin{prop}\label{onebound}
    For a $G$-Mackey functor $\uM$, if $\uM(G/e)=0$, then $\Sigma^{1} H\uM \in [2, |G|].$
\end{prop}

\begin{proof}
   From the previous proposition, we have $\Sigma^{1} H\uM \ge 2.$ For the other bound,  we need to prove that the group  $\tH^{1-k-r}_H(S^{k\bar{\rho}_H}; \uM)$ vanishes for each subgroup $H \le G$ and $k|H|> |G|$ and $r\ge 0.$ This forces $k \ge 2.$ Thus, for any $r\ge 0$, $1-k-r<0$, hence the group vanishes for all $H \le G.$ Thus the result follows.
\end{proof}

To obtain an alternative characterization of $2$-slices of the form
$\Sigma^1 H\uM$, we begin with the following definitions.

\begin{defn}
  For a finite group $G$ and $k \ge 1$, define $\cP^{(k)}_G:=\{L \in \cP_G: |L|=k\}$. For a Mackey functor $\uM$, we define a map $\psi_H: \uM(G/H)\to \displaystyle\prod_{C\in \cP^{(2)}_H}\uM(G/C)$ given by $\psi_H(m)=(\res^H_C(m))_C$.
\end{defn}
With this, we can completely derive the algebraic criterion for $\Sigma^{1} H\uM$ to be a $2$-slice as follows.
\begin{prop}\label{fullsigman}
   $\Sigma^1H\uM$ is a $2$-slice if and only if $\psi_H$ is injective for all $H \le G.$
\end{prop}
\begin{proof}
    By \eqref{kervarphi} and Proposition \ref{onesus}, it is enough to prove that the following are equivalent.
    \begin{enumerate}
        \item $\uM(G/e)=0$ and $\bigcap_{L \in \cP_H}\ker\res^H_L=0$ for $H \le G$ with $|H|\ge 3.$
        \item $\psi_H$ is injective for all $H \le G.$
    \end{enumerate}
   Assume that $(2)$ holds. Then for $H=e$, we have $\psi_e: \uM(G/e)\to 0$ is injective. Hence $\uM(G/e)=0.$
   Assume $|H| \ge 3$, then $\bigcap_{L \in \cP_H}\ker\res^H_L \subseteq \ker \psi_H$. Thus (1) holds.

   On the other hand, assume that (1) is true. To prove $\psi_H$ is injective for all $H\le G$, we use induction on $|H|.$ Note that $|H|=1, 2$ are obvious. Suppose $|H| \ge 3.$ By induction, we assume that $\psi_L$ is injective for all $L \in \cP_H.$ Take $m \in \ker\psi_H$. Then $\psi_H(m)= (\res^H_L(m))_{L}=0$, where $L \in \cP_H^{(2)}.$ Then $\psi_L(\res^H_L(m))$ are precisely the components of
$\psi_H(m)$. Hence, $\psi_L(\res^H_L(m))=0$, and by induction,
we obtain $\res^H_L(m)=0$. Therefore, $m\in\bigcap_{L\in\cP_H}\ker(\res^H_L)=0.$ Hence $\psi_{H}$ is injective.
\end{proof}

The following is an immediate corollary.

\begin{cor}\label{oddvanish}
    If $\Sigma^1H\uM$ is a $2$-slice, then $\uM(G/H)=0$ for every subgroup
$H\le G$ of odd order.
\end{cor}

\begin{proof}
    Let $H$ be a subgroup of $G$ of odd order. By Lagrange's theorem there is no subgroup $C$ of $H$ with $|C|=2$, hence $\cP_H^{(2)}=\emptyset.$ This proves the result.
\end{proof}

\begin{remark}
    The converse of Corollary \ref{oddvanish} fails. For example, let $G$ be any group with even order $|G|\ge 4.$ Consider $\uM=\ug$, given by $\ug(G/G)=\F$ and $\ug(G/L)=0$ for $L\in \cP_G.$ Then $\uM(G/K)=0$ for all odd order subgroup $K$ of $G.$ However $\Sigma^{1} H\uM$ is not a $2$-slice because
    \[
    [S^{\rho_G}, \Sigma^{1} H\uM]^G\cong \tH^{0}_G(S^{\bar{\rho}_G}; \uM)\cong \ker(\uM(G/G) \to \tH^0_G(S(\bar{\rho}_G)_+; \uM))\cong \F\ne 0.
    \]
    The last equality follows from Theorem \ref{zerosrho}.
\end{remark}

Having introduced the necessary tools, we are now ready to state our main theorem, which characterizes \(2\)-slices for an arbitrary finite group \(G\).
\begin{thm}\label{main2slice}
    A $G$-spectrum $Y$ is a $2$-slice if and only if the following conditions hold.
    \begin{enumerate}
        \item $\upi_k Y=0$ for $k \neq 1, 2$ and $\upi_1 Y(G/e)=0.$

        \item $\Sigma^2 H \upi_2(Y)$ is a $2$-slice.

        \item For every $H \le G$ with $|H|\ge 3,$ the $k$-invariant $\kappa: \Sigma^1 H\upi_1 Y \to  \Sigma^3 H\upi_2 Y$ induces an injective map $\ker(\varphi_H(\upi_1 Y )) \to \tH^1_H(S(\bar{\rho}_H)_+; \upi_2Y).$
    \end{enumerate}
\end{thm}

\begin{proof} For the proof, we first assume that $Y$ is a $2$-slice.

    Since \(Y\geq 2\), Theorem~\ref{lbound} implies that $
\pi_k(Y^e)=0$ for $k<2$.
Consequently, \(\upi_1Y(G/e)=0\) and \(\upi_kY=0\) for \(k\leq 0\). While for  \(k\geq 3\), we have \(S^k\geq 3\). Thus, $[S^k,Y]^G=0$. 
It follows that the only possibly non-vanishing homotopy Mackey functors of \(Y\) are \({\upi}_1Y\) and \({\upi}_2Y\).
    Hence, \(Y\) fits into the following fiber sequence:
\begin{myeq}\label{2-cofiber}
\Sigma^2 H\upi_2Y
\longrightarrow Y
\longrightarrow \Sigma^1 H\upi_1Y.
\end{myeq}

    Since $H \upi_1Y \le 0$ and $Y \le 2$, therefore we have  $\Sigma ^2H \upi_2Y\le 2.$ On the other hand $\Sigma ^2H \upi_2Y \ge 2$ always holds. Hence $\Sigma ^2H \upi_2Y$ is a $2$-slice.

    Consider the long exact sequence induced from the above fiber sequence
    \[
    \xymatrix@C=.6cm{\cdots  [S^{k\rho_H+r}, \Sigma^2 H \upi_2Y]^H \ar[d]^{\cong} \ar[r] & [S^{k\rho_H+r}, Y]^H \ar[r] & [S^{k\rho_H+r}, \Sigma H \upi_1Y]^H \ar[d]^{\cong } \ar[r] & [S^{k\rho_H+r}, \Sigma^3 H \upi_2Y]^H \ar[d]^{\cong } \cdots \\\tH^{2-r-k}_H(S^{k\bar{\rho}_H}; \upi_2Y) & & \tH^{1-r-k}_H(S^{k\bar{\rho}_H}; \upi_1Y) \ar[r] & \tH^{3-r-k}_H(S^{k\bar{\rho}_H}; \upi_2Y)}
    \]
    Let $k |H|>2$ and $r \ge 0.$ Hence the group $\tH^{1-r-k}_H(S^{k\bar{\rho}_H}; \upi_1Y)$ could be non-zero only for the admissible pair $(k,r)=(1,0)$, and then this group is $\ker(\varphi_H(\upi_1 Y ))$ and the other group $\tH^{2}_H(S^{\bar{\rho}_H}; \upi_2Y)\cong \tH^1_H(S(\bar{\rho}_H)_+; \upi_2Y)$ and the left corner vanishes since $\Sigma ^2H \upi_2Y\le 2$. Hence $[S^{k\rho_H+r}, Y]^H \cong \ker(\varphi_H(\upi_1 Y )) \to \tH^1_H(S(\bar{\rho}_H)_+; \upi_2Y).$

    For the converse part, let us assume that all these conditions hold. The condition (3) suggests that $Y\le 2.$ Again consider the fiber sequence~\eqref{2-cofiber}, since by condition (1), we have $\upi_1Y(G/e)$. Therefore, by Proposition~\ref{onebound}, we obtain
$\Sigma^1 H\upi_1Y\geq 2$.

    By condition (2),  $\Sigma^2 H \upi_2(Y) \ge 2.$ Since the localizing subcategory is closed under extensions hence $Y \ge 2.$ 
\end{proof}

\begin{mysubsection}{Alternating group $A_4$}~ Let \(G=A_4\) be the alternating group of order \(12\). Up to conjugacy, the subgroups of \(A_4\) are the trivial subgroup \(e\), a cyclic subgroup \(C_2\) of order \(2\), a cyclic subgroup \(C_3\) of order \(3\), the Klein four-group \(K_4\), which is normal in \(A_4\), and \(A_4\) itself.  Assume that $H_1,H_2,H_3$ are the three subgroups of order $2$ in $K_4$ (and hence in $A_4$). Since $A_4$ has no
subgroup of index $2$, the groups $A_i$ in
Corollary~\ref{zerosrho-max} are simply Weyl fixed points.  



We first consider the Eilenberg–Mac Lane spectra. Applying Propositions~\ref{twosus} and \ref{onesus} to the group $A_4$, we obtain the following results.
\begin{prop}
Let $\uM$ be an $A_4$-Mackey functor. 
\begin{enumerate}
   \item The spectrum \(\Sigma^2 H\uM\) is a \(2\)-slice if and only if the following conditions are satisfied:
\begin{enumerate}
    \item The restriction map $\res^{C_2}_e:\uM(A_4/C_2)\longrightarrow \uM(A_4/e)$ is injective.

    \item The restriction map $\res^{C_3}_e:\uM(A_4/C_3)\longrightarrow \uM(A_4/e)^{C_3}$    is bijective.

    \item The map $\bigl(\res^{K_4}_{H_1},\res^{K_4}_{H_2},
    \res^{K_4}_{H_3}\bigr):
    \uM(A_4/K_4)\longrightarrow T_{K_4}(\uM)$
    is an bijective.

    \item The map $\bigl(\res^{A_4}_{K_4},\res^{A_4}_{C_3}\bigr):
    \uM(A_4/A_4)\longrightarrow T_{A_4}(\uM)$
    is bijective.
\end{enumerate}
Here, $T_{K_4}(\uM):=
\tH^0_{K_4}\bigl(S(\rb_{K_4})_+;\uM\bigr)$
and $T_{A_4}(\uM):=
\tH^0_{A_4}\bigl(S(\rb_{A_4})_+;\uM\bigr),$
where these groups are described in Corollaries~\ref{algk4} and~\ref{alga4}, respectively.
\item The spectrum $\Sigma^1 H\uM$ is a $2$-slice if and only if the following conditions are satisfied:
    \begin{enumerate}
        \item The Mackey functor value at orbits $A_4/e$ and $A_4/C_3$ vanish. ie.
        \[
        \uM(A_4/e)=0 \text{ and } \uM(A_4/C_3)=0.
        \]
        \item The map $\uM(A_4/K_4) \to \displaystyle\bigoplus_{i=1}^3 \uM(A_4/H_i)$ is injective.
        \item The restriction map $\res^{A_4}_{K_4}$ is injective.
\end{enumerate}
\end{enumerate}
\end{prop}
\begin{proof} (1) This result follows directly from Proposition~\ref{twosus}, together with Corollaries~\ref{algcpn} and~\ref{alga4}.

(2) This assertion follows from Proposition~\ref{onesus}.
\end{proof}

For a general $A_4$-spectrum $Y$, the full classification is as follows.

\begin{thm}
    An $A_4$-spectrum $Y$ is a $2$-slice if and only if 
\begin{enumerate}
    \item The homotopy Mackey functor $\upi_{k}Y$ are zero for $k\neq 1,2$ and  $\upi_1Y(A_4/e)=0$.
    \item The spectrum $\Sigma^2 H\upi_2 Y$ is a 2-slice. 
    \item The $\kappa$-invariant $\Sigma^1 H\upi_1 Y \to \Sigma^3 H\upi_2 Y$ induces the following injections:
    \begin{align*}
\upi_1Y(A_4/C_3)&\hookrightarrow H^1_{C_3}\bigl(S(\bar\rho_{C_3});\res_{C_3}\upi_2Y\bigr)=\upi_2Y(A_4/e)_{C_3},\\
\bigcap_{i=1}^3\ker(\res^{K_4}_{H_i})&\hookrightarrow
\tH^1_{K_4}\bigl(S(\bar\rho_{K_4})_+;\res_{K_4} \upi_2 Y\bigr),\\
\ker(\res^{A_4}_{K_4})\cap\ker(\res^{A_4}_{C_3})&\hookrightarrow
\tH^1_{A_4}\bigl(S(\bar\rho_{A_4})_+; \upi_2 Y\bigr).
\end{align*}
\end{enumerate}
\end{thm}
\begin{proof}
    This characterization is the direct consequence of Theorem~\ref{main2slice} along with Corollaries~\ref{algcpn}, \ref{algk4} and~\ref{alga4}.
\end{proof}
\end{mysubsection}

\begin{mysubsection}{Cyclic $p$-groups}~ For a prime $p$, consider the cyclic $p$-group $C_{p^n}$.
It has a unique maximal subgroup, namely $C_{p^{n-1}}$.
Since the subgroup lattice of $C_{p^n}$ is a chain, each subgroup contributes at most one condition. Thus, we can index these conditions by $j$.

Applying Propositions~\ref{twosus} and \ref{onesus} to the group $C_{p^n}$, we obtain the following results.

\begin{prop}
Let $\uM$ be a $C_{p^n}$-Mackey functor. Then the following hold.

\begin{enumerate}[label=\textup{(\arabic*)}]

\item Suppose that $p$ is odd. Then $\Sigma^2H\uM$ is a $2$-slice
if and only if, for every $1\leq j\leq n$, the restriction map
\[
\res^{C_{p^j}}_{C_{p^{j-1}}}\colon
\uM(C_{p^n}/C_{p^j})
\longrightarrow
\uM(C_{p^n}/C_{p^{j-1}})^{C_{p^j}/C_{p^{j-1}}}
\]
is bijective.

\item Suppose that $p=2$. Then $\Sigma^2H\uM$ is a $2$-slice
if and only if the following conditions hold:
\begin{enumerate}[label=\textup{(\alph*)}]
    \item If $n=1$, then $\res^{C_2}_e$ is injective.
    \item If $n\geq 2$, then $\res^{C_2}_e$ is injective, and for
    every $2\leq j\leq n$, the restriction map $\res^{C_{2^j}}_{C_{2^{j-1}}}$
    is a bijection to the following set
    \[
    \left\{x\in\uM(C_{2^j}/C_{2^{j-1}})\ \middle|\
    \res^{C_{2^{j-1}}}_{C_{2^{j-2}}}(x)
    \in
    \uM(C_{2^j}/C_{2^{j-2}})^{C_{2^{j}}/C_{2^{j-2}}}
    \right\}.
    \]
\end{enumerate}

\item Suppose that $p=2$. Then $\Sigma^1 H\uM$ is a $2$-slice
if and only if the following conditions hold:
\begin{enumerate}[label=\textup{(\alph*)}]
    \item If $n=1$, then $\uM(C_2/e)=0$.
    \item If $n\geq 2$, then $\uM(C_{2^n}/e)=0$, and for every
    $2\leq j\leq n$, the restriction map $\res^{C_{2^j}}_{C_{2^{j-1}}}$
    is injective to the following set
    \[
    \left\{x\in\uM(C_{2^j}/C_{2^{j-1}})\ \middle|\
    \res^{C_{2^{j-1}}}_{C_{2^{j-2}}}(x)
    \in
    \uM(C_{2^j}/C_{2^{j-2}})^{C_{2^j}/C_{2^{j-2}}}
    \right\}.
    \]
\end{enumerate}

\end{enumerate}
\end{prop}
\begin{proof}
    The result follows directly from Corollary~\ref{algcpn} and Proposition~\ref{twosus}.
\end{proof}

As an immediate consequence of Proposition~\ref{fullsigman}, we obtain the following remark for cyclic groups of odd prime-power order.

\begin{remark}
Let $p$ be an odd prime and let $\uM$ be a nonzero $C_{p^n}$-Mackey functor. Then $\Sigma^1 H\uM$ cannot be a $2$-slice.
\end{remark}

For a $C_{p^n}$-spectrum
$Y$, write
\[
\kappa_j\colon
\ker\bigl(\res^{C_{p^j}}_{C_{p^{j-1}}}\colon
\upi_1 Y(C_{p^n}/C_{p^j})\to \upi_1 Y(C_{p^n}/C_{p^{j-1}})\bigr)\longrightarrow
\tH^1_{C_{p^j}}\bigl(S(\bar\rho_{C_{p^j}})_+;\res_{C_{p^j}}\upi_2 Y\bigr)
\]
  induced by the $\kappa$-invariant $\kappa: \Sigma H\upi_1 Y \to \Sigma^3 H\upi_2Y$.

  Applying Theorem~\ref{main2slice} to the cyclic $p$-group $C_{p^n}$, we obtain the following characterization of $2$-slices.   
\begin{thm}
    A $C_{p^n}$-spectrum $Y$ is a $2$-slice if and only if 
\begin{enumerate}
    \item $\upi_k Y\cong 0$ if $k\neq 1,2$ and $\upi_1Y(C_{p^n}/e)=0$.
    \item $\Sigma^2 H\upi_2 Y$ is a 2-slice.
    
    \item If $p$ is odd, for every $1\le j\le n$, $\kappa_j$ is injective.

    \item If $p=2$ and $n\geq 2$, then $\kappa_j$ is injective for every $2\leq j\leq n$.
\end{enumerate}
\end{thm}
\end{mysubsection}

\begin{remark}
The computation of $\tH^1_{G}(S(\bar{\rho}_G)_+; \uM)$ can be done following the paper \cite[(3.5)]{HHR16}. More precisely, let $S_k$ denote the flag $F_0 \subsetneq \cdots \subsetneq F_k$ of $G$-sets and $G$ acts by left translation. We have a face map $d_i: S_k \to S_{k-1}$ that drops $F_i$. For a $G$-Mackey functor, this induces a cochain complex 
\[
\uM(G/G)\stackrel{\delta^0}{\to} \uM(S_0)\stackrel{\delta^1}{\to}\uM(S_1) \stackrel{\delta^2}{\to}\uM(S_2) \cdots.
\]
Cohomology of this cochain complex provides $\tH^k_{G}(S(\bar{\rho}_G)_+; \uM)$ via the isomorphism 
\[
\tH^k_{G}(S(\bar{\rho}_G)_+; \uM) \cong \tH^{k+1}_{G}(S^{\bar{\rho}_G}; \uM) \text{ for all } k\ge 1.
\]
One may compute the same using the filtration of the Mackey functor $\uM$ following Grodal's work \cite{Gro02}.
\end{remark}


\begin{thebibliography}{99}

\bibitem{Dug05}
{\scshape D.~Dugger,} \emph{An {A}tiyah--{H}irzebruch spectral sequence for
{$KR$}-theory}, K-Theory \textbf{35} (2005), no.~3-4, 213--256.

\bibitem{Gro02}
{\scshape J.~Grodal,} \emph{Higher Limits via Subgroup Complexes},  Ann. \ of Math. \ (2) Second Series, {\bf 155}, No. 2 (2002), pp. 405-457

\bibitem{HHR16}
{\scshape M.~A.~Hill, M.~J.~Hopkins, D.~C.~Ravenel,}
\emph{On the nonexistence of elements of Kervaire invariant one,}
Ann.\ of Math.\ (2) {\bf 184} (2016), no.~1, 1--262.


\bibitem{HY18}
{\scshape M.~A.~Hill, C.~Yarnall,}
\emph{A new formulation of the equivariant slice filtration with applications to $C_p$-slices,}
Proc.\ Amer.\ Math.\ Soc.\ {\bf 146} (2018), no.~8, 3605--3614.

\bibitem{MNN19}
{\scshape A. ~Mathew, N. ~Naumann, Niko,  J.~Noel, Justin,}
\emph{Derived induction and restriction theory,} Geom. Topol. {\bf 23} (2019), no. 2, 541–636.

\bibitem{Slo22}
{\scshape C.~Slone,}
\emph{Klein four $2$-slices and the slices of $\Sigma^{\pm n}H\underline{\mathbb Z}$,}
Math.\ Z.\ {\bf 301} (2022), 3895--3938.

\bibitem{Voe02}
{\scshape V.~Voevodsky,} \emph{Open problems in the motivic stable homotopy
theory.~{I}}, Motives, polylogarithms and {H}odge theory, {P}art {I}
({I}rvine, {CA}, 1998), Int. Press Lect. Ser., vol.~3, Int. Press,
Somerville, MA, 2002, pp.~3--34.

\bibitem{U}
{\scshape J.~R.~Ullman,}
\emph{On the regular slice spectral sequence,}
Ph.D.\ thesis, Massachusetts Institute of Technology, 2013.



\end{thebibliography}
\end{document}